# DISCUSSION: CONDITIONAL GROWTH CHARTS

By Anneli Pere

*University of Helsinki*

First, I would like to congratulate the authors for an interesting application of their semiparametric quantile regression model to longitudinal human growth data. On an earlier occasion, I had an opportunity to collaborate with the authors on applying this method for the purpose of constructing semilongitudinal growth charts for height and for body mass index. As my background is in medicine rather than statistics, I will focus my comments below more on the underlying biological aspects and leave the mathematical and statistical comments to the other discussants.

**1. Human growth.** Much of the current understanding of physical growth of children derives from auxological works from as early as the 1950s. They are still quite valid because the biological nature of human growth has remained basically the same. Two good sources for understanding human growth can be pointed out: the book on the history of ideas concerning growth written by James Tanner [4] and another book with general text, written by David Sinclair [2], that gives an overview of the various manifestations of human growth.

Human growth can be divided into the phases of fetal, infant, childhood and pubertal growth. These phases overlap in time and interact with each other, that is, development during one phase may influence that in another. A child inherits, separately from both parents, the genes that largely determine the "growth channel," the tempo or timing of growth events and the potential for adult height. Potential adult height can be, to a good approximation, viewed as a built-in or "programmed" property. If the growth of a child has been disturbed for a period that is not exceedingly long, by an environmental factor or an illness, some form of catch-up or catch-down growth usually follows. For this reason, it seems unlikely that final adult height can be increased much, for example, by medication. As an example of an inherited growth pattern, a child can be long and stout during infancy but gradually become slender and shorter in stature in early childhood. Another example is a child whose one parent was late in maturation during puberty and is tall as an adult, whereas the other parent was early in maturation and is short. This child could have inherited, say, late pubertal maturation and short adult height and would be exceptionally short at the age when most of his/her peers have entered puberty and their growth has accelerated according to their pubertal growth spurt. Assessment of growth during puberty is difficult without any knowledge about the "biological" or "maturational" age of the child. All in all, growth is a complicated process with a series of changes, not just addition of material.

**2. Growth surveillance in Finland.** In Finland, trained nurses have carried out follow-up of height and weight growth of practically all children in the well-baby clinics and in the schools now for almost half a century. Growth curves of the parents of the children in follow-up today are already on record and available for researchers, and soon also those of their grandparents will be. Growth curves of siblings can also be assessed for possible peculiarities in the growth patterns.

A new form of growth charts was presented twenty years ago [3] to facilitate early detection of aberrant growth. Height in these charts is presented as a standard deviation score (SDS), that is, deviation of height in SD units from the mean height for each considered age and sex. Weight correlates more strongly with height than with age. Therefore weight against height, rather than weight against age, is used in these charts, and it is presented as a deviation, in percentage, from the median weight for the considered height and sex. Knowledge of weight-for-height, say 20% below the median, and of age, suggests to an experienced clinician immediately an image of the body build of the child in question: an underweight child! Therefore individual changes in the body build are easily detected from the weight-for-height curve. On the other hand, if weight-for-age is used, interpretations of body build will generally depend on whether the child is short, average or tall. In the Finnish growth charts, an approximately horizontal line represents normal growth. Any changes from an individual's "growth channel," such as an upward or a downward bend, can easily be detected by eye (compare the traditional and the new format for height chart in Figure 1). Moreover, magnitudes of such changes are then comparable across different ages, unlike in a situation in which the original measurements (in centimeters/inches or kilograms/pounds) are retained. The weight-for-height percentage values can be plotted on the height chart, to the corresponding measurement ages, but the $y$-scale of height SDS needs to be multiplied by 10 to fit better the weight-for-height scale. From this chart the simultaneous growth of height and weight can be readily assessed (see Figures 2 and 3).

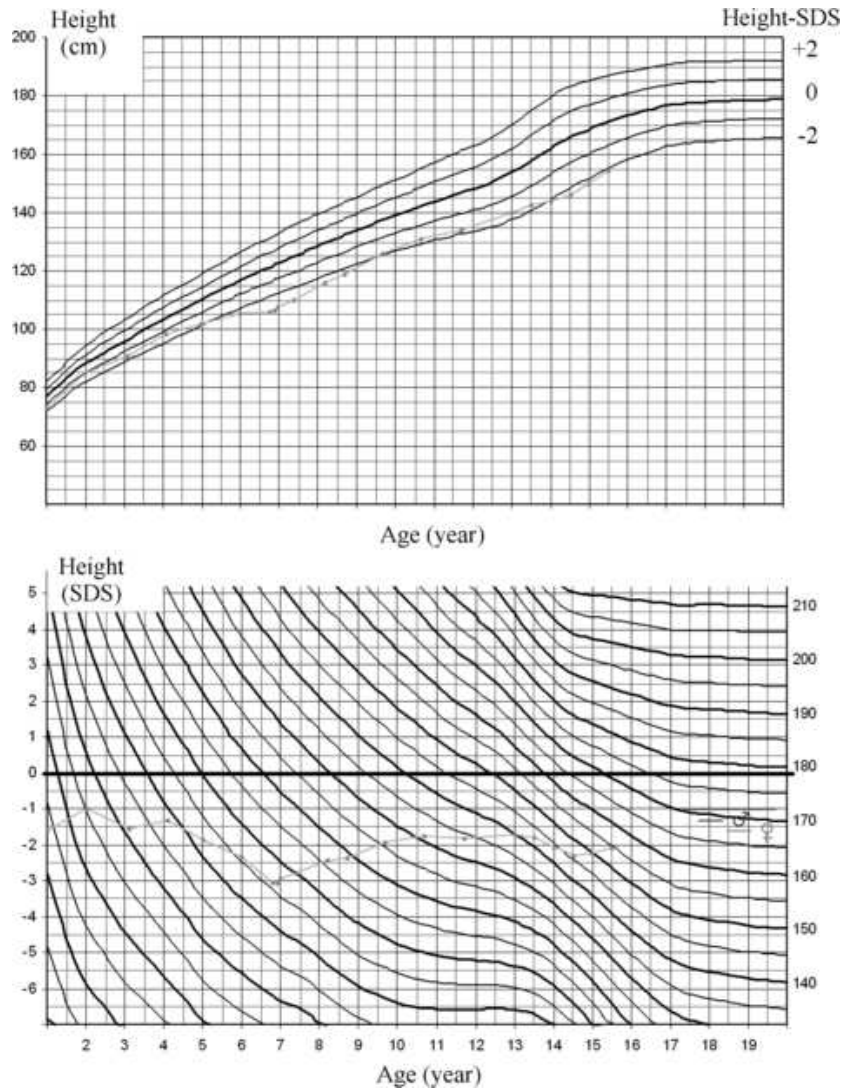

Fig. 1. *An example of pathological growth plotted in old-fashioned (top panel) and in current format (bottom panel) of the Finnish height chart. On the current height charts, the horizontal lines and y-axis scale represent height in SDS (standard deviation score) and the curved lines in the background indicate the absolute heights (in cm). This boy was diagnosed to have hypothyroidism (lack of thyroid hormone) caused by autoimmune thyroiditis at the age of 6.9 years. A clear bend in height growth had appeared at 4.1 years, more than two years before other clinical signs typical for this disease had appeared. The growth aberration is not as easily detected from the old-fashioned chart. Catch-up growth appeared soon after the medication (substitution with thyroid hormone) was started. Note the downward bend in height is away from the expected height-SDS, −1.05 SDS, calculated from the parents' heights (mother's −1.65 SDS and father's −1.4 SDS). The short horizontal lines on the right indicate these three adult height-SDS's.*

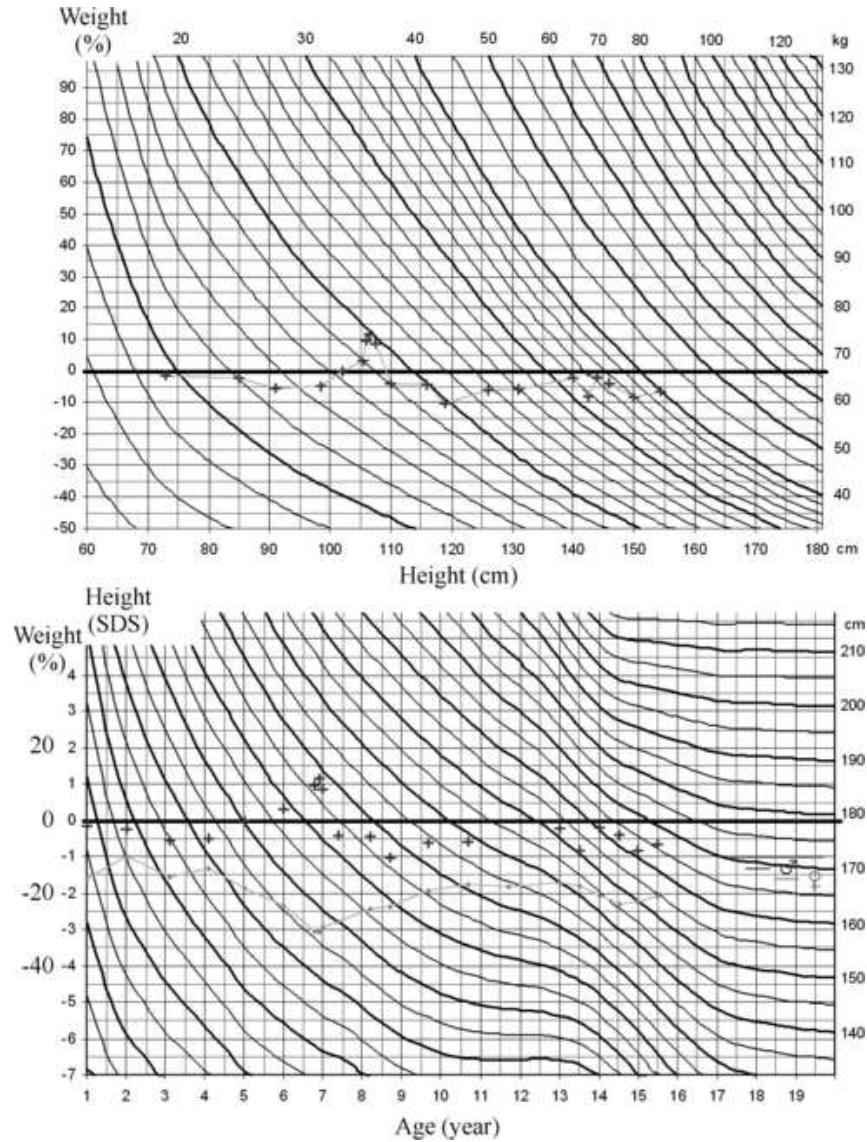

FIG. 2. *On the current Finnish weight-for-height chart (top panel), the horizontal lines and y-axis scale indicate weight presented as the percentage deviation from the median weight for the specific height of the same sex. The absolute weight (in kg) can be read from the curved background lines. An upward bend in the weight curve can be observed; the body build changes from a bit underweight (weight-for-height 5% below the median, 14.9 kg at 98.5 cm height) to a bit plump (12% above the median, 19.8 kg at 106.5 cm height). The ages when the measurements were done cannot be read from this chart. For this purpose, the weight-for-height measurements can be plotted on the height-for-age chart as shown in the bottom panel: crosses indicate weight-for-height and dots, connected with lines, indicate height-SDS. The change in body build would have been missed if a weight-for-age chart had been used: the increase in absolute weight was only 4.9 kg during 2.8 years.*

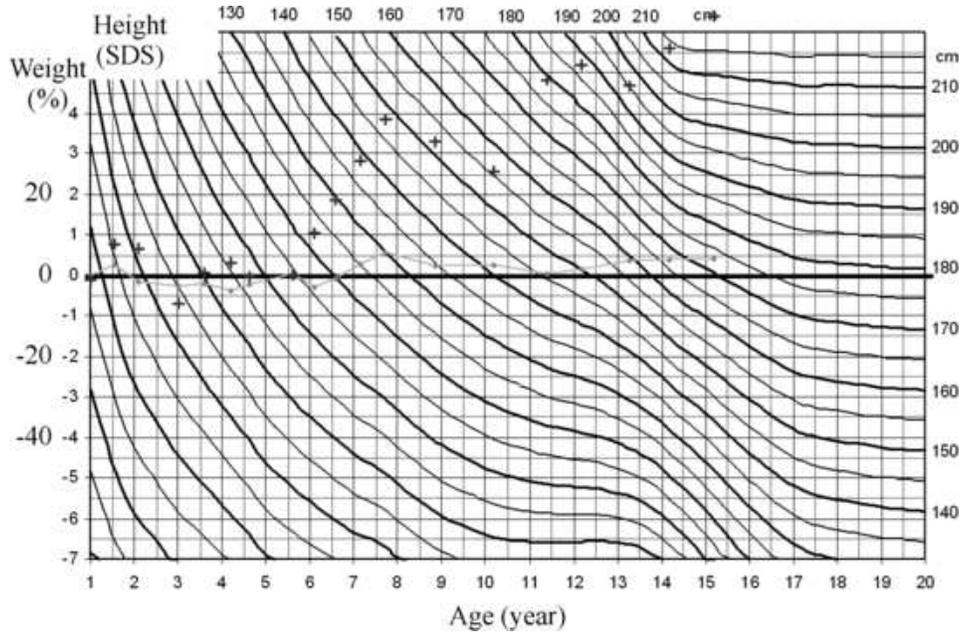

FIG. 3. *This example of growth of an obese but otherwise healthy boy shows how difficult it is to distinguish exceptional but healthy growth from pathological growth with an underlying disease. The percentage values for weight-for-height have been plotted to this height chart in the same manner as described for the bottom panel in Figure* 2. *This makes simultaneous assessment of height* (*dots connected with lines*) *and weight* (*crosses*) *growth easier: in this example, soon after the steep upward bend in the weight curve* (*at* 5.61 *years*) *height growth accelerates a bit* (*height curve "bending" at* 6.10 *years*). *This combination of bends, and particularly when the bend in the weight curve is not continuous, is a less alarming pattern of growth deviation than the one presented in Figure* 2.

In growth surveillance, the use of population averages is not that important. The fundamental idea is to understand what is normal for a particular child. Clear changes from the previous growth path are alarming and they should be noticed before leading to exceptional size. New guidelines for screening for pathological growth were given ten years ago. The tables for screening limits are printed on the growth charts and also included in the software (developed by Markkula [1] in 1996–2005) used in many centers to register growth measurements and to plot individual growth curves (such as those in Figures 1–3). These screening limits are based on estimates of the extreme 0.5% of healthy children who increase or decrease their height SDS or weight-for-height percentage units. In other words, these limits are based on children whose growth of height or weight has been exceptionally fast or slow. At any given screening age, 99% of the healthy children are considered to stay within their height SDS path, and likewise 99% to stay within their individual weight-for-height path. A potential problem in screening is that it is not always the same individuals who are among the 99% that are considered healthy. Therefore, the more screenings are carried out, the larger proportion of healthy children are classified as exceptional. On the other hand, there will also always be false negatives, that is, not all children with a disease that causes pathological growth can be screened with these screening limits. There also remains much in the recognition of patterns that cannot be put to figures. Some guidelines have been added to the tables for screening limits in words. For example, a change in weight-for-height toward the median of the population, or a change in height-SDS toward the mean of the population or toward the expected height SDS calculated from the parents' heights is less alarming than changes in the opposite direction (see Figure 1). Likewise, a slow but continuous change in the individual path (a steady, nonwavy bend in the growth curve) is more alarming than an unsteady change which already shows some correction of the disturbance (see Figure 3). Could these two guidelines be implemented to the conditional quantile regression model?

**3. The global quantile regression model and the example of infant weight.** When considering complex developments such as human growth, it is difficult to give to the global model considered in (2.1) an interpretation that would make biological sense. As I have outlined above, one can hardly justify an application of a simple linear model directly based on past measurements of height or weight. Moreover, when considering one of these measurements—in the example weight—it seems rather odd to simply add the child's height as is done in (5.1), albeit multiplied with some coefficient, to such a linear predictor. Follow-up of length/height growth has more value than follow-up of weight growth when it comes to early detection of certain diseases (that may cause a growth aberration before signs typical for the disease appear). Clearly, given the target of detecting aberrant growth, it would be desirable that much more information from the child's past growth history would be included in the predictor than is done in the example. In particular, it would seem to be important to include the previous growth pattern such as steady or unsteady change in the individual's quantile position.

On the other hand, if predictions provided by quantile regression are expected to discern between normal and abnormal growth, they should not be "too good." Otherwise it may well happen that if abnormal growth of a child develops slowly, that is, insidiously, and is not detected by the method at an early stage, the consequent individual predictions will say only that "the same type of growth is likely to continue." Such predictions are nearly self-fulfilling and pathological growth is not distinguished from unusual but healthy growth. From this perspective, it may be best if the model is estimated using only growth data for healthy children who serve as a good reference population. Growth data for the cases of some rare, treatable diseases wished to be screened at an early stage are always collected retrospectively and are also bound to contain measurement or recording errors and such.

In the particular example considered, development of the weight of a boy was considered from birth to an age of 0.61 year. This boy seemed to have postnatal catch-down growth of weight during the first two months, after which he grew steadily a little below the median. The measurement at 0.46 year is somewhat below the earlier growth path, and the next one, at 0.61 year, is above the path. Overall, however, this early infancy growth pattern does not seem unusual because more than two consecutive measurements should be assessed. The downward bend away from the median at the age of 0.46 year (if it is not a measurement error, which is also a real possibility!) is actually clinically a more alarming sign than the upward bend thereafter. It is highly unlikely to have clinically meaningful "bends" in opposite directions and so close to each other. Rather, if not measurement errors, these readings only reflect natural variation in the speed of weight growth in a healthy infant. In fact, for a clinically meaningful conclusion, one would have to assess growth in length and weight simultaneously.

Two phenomena can be detected from the measurement at 0.61 year: regression toward the mean (of the population) and the child's own "tracking" or returning to his growth path. Most likely, the next weighing of this child would have already shown deceleration of weight growth. Regression toward the mean can be nicely seen in the comparison between the global model and the LMS-AR model (Figure 2 in the paper). The LMS-AR model predicts that the child's individual growth quantile remains the same and that deviations toward the population mean or away from it are exceptional. Contrary to this, the global model shows regression toward the mean: even the 0.5 quantile is closer to the population average weight than the respective quantile with the LMS-AR model, and the difference between the two methods becomes greater the closer the population average (approximately 9.0 kg, estimated from Figure 1 of the paper) is.

**4. Future work in screening for pathological growth.** Information on both conditional height growth and conditional weight growth is needed to give a picture of how exceptionally unsteadily growing but healthy children grow. Similar studies should then be performed with growth data on children with certain diseases that could be diagnosed earlier if proper growth surveillance were carried out. Would the results of an application of the global model be easier to interpret if the quantile position rather than age and the absolute value (of height or weight) were used?

Most interesting, from a biological and also a clinical perspective, would be analyses based on a joint consideration (two-dimensional distribution) of height and weight at various ages. I am thinking of a "topographical map" where the median quantile is in the middle and other quantiles form contours of various shapes around the middle. It would be informative to learn how healthy children change their position with age on this map.

**5. Other application areas in medicine.** Many laboratory reference values could be developed with quantile regression models that do not make assumptions about the distribution. The conditional model is especially interesting and could be applied to common laboratory values such as blood hemoglobin that are known to have individual levels. For the child population, several laboratory reference values are age-specific and here, too, the use of the conditional model could be fruitful. One practical problem remains, however: the size of the data set used for many reference values, especially for children, is often too small for proper estimation of the extreme quantiles.

**Acknowledgment.** The author thanks Professor Elja Arjas for valuable help in preparing this discussion.

THE HOSPITAL FOR CHILDREN AND ADOLESCENTS
UNIVERSITY OF HELSINKI
STENBÄCKINKATU 9
FI-00029 HELSINKI
FINLAND
E-MAIL: anneli.pere@fimnet.fi